\documentclass[CJK]{article}
\usepackage[dvips]{color}
\usepackage{amsmath}
\usepackage{amsfonts}   
\usepackage{amssymb}    
\usepackage{multicol}
\usepackage{mathptm,graphicx}
\usepackage{epsf}
\newtheorem{theorem}{Theorem}[section]
\newtheorem{lemma}{Lemma}[section]
\newtheorem{remark}{Remark}[section]
\newtheorem{definition}{Definition}[section]

\newtheorem{proposition}{Proposition}[section]
\newtheorem{assum}{Assumption}[section]

\newcommand{\sce}[1]{\setcounter{equation}{#1}}

\date{}
\begin{document}
\title{G-Doob-Meyer Decomposition and Its Application
in Bid-Ask Pricing for American Contingent Claim Under Knightian
Uncertainty
}
\author{Wei Chen \\
Institute of Quantitative Economics\\
 School of Economics\\
    Shandong University\\
    250199, Jinan, China\\
weichen@sdu.edu.cn\\
}
\date{}
\maketitle
\begin{center}
\begin{minipage}{120mm}
\baselineskip 0.2in {\small {\bf Abstract} The target of this
paper is to establish the bid-ask pricing frame work for the
American contingent claims against risky assets with G-asset price
systems (see \cite{Chen2013b}) on the financial market under
Knight uncertainty. First, we prove G-Dooby-Meyer decomposition
for G-supermartingale. Furthermore, we consider bid-ask pricing
American contingent claims under Knight uncertain, by using
G-Dooby-Meyer decomposition, we construct dynamic superhedge
stragies for the optimal stopping problem, and prove that the
value functions of the optimal stopping problems are the bid and
ask prices of the American contingent claims under Knight
uncertain. Finally, we consider a free boundary problem, prove the
strong solution existence of the free boundary problem, and derive
that the value function of the optimal stopping problem is
equivalent to the strong solution to the free boundary problem.
}
\end{minipage}
\end{center}
\newpage
\section{Introduction}

The earliest, and one of the most penetrating, analysis on the
pricing of the American option is by McKean \cite{McKean}. There
the problem of pricing the American option is transformed into a
Stefan or free boundary problem. Solving the latter, McKean writes
the American option price explicitly up to knowing a certain
function -- the optimal stopping boundary.

Bensoussan \cite{Bensoussan} presents a rigorous treatment for
American contingent claims, that can be exercised at any time
before or at maturity. He adapts the Black and Scholes
\cite{Black} methodology of duplicating the cash flow from such a
claim to this situation by skillfully managing a self-financing
portfolio that contains only the basic instruments of the market,
i.e., the stocks and the bond, and that entails no arbitrage
opportunities before exercise. Bensoussan shows that the pricing
of such claims is indeed possible and characterized the exercise
time by means of an appropriate optimal stopping problem. In the
study of the latter, Bensoussan employs the so-called
"penalization method", which forces rather stringent boundedness
and regularity conditions on the payoff from the contingent claim.

From the theory of optimal stopping, it is well known that the
value process of the optimal stopping problem can be characterized
as the smallest supermartingale majorant to the stopping reward.
Base on the Doob-Meyer decomposition for the supermartingale, a
"martingale" treatment of the optimal stopping problem is used for
handling pricing the American option by Karatzas \cite{Karatzas},
EL Karoui and Karatzas \cite{KarouiKaratzas1989},
\cite{KarouiKaratzas1991}.

The Doob decomposition Theorem was proved by and is named for
Joseph L. Doob \cite{Doob1953}. The analogous theorem in the
continuous time case is the Doob-Meyer decomposition theorem
proved by Meyer in \cite{Meyer1962} and \cite{Meyer1963}. For the
pricing American option problem in incomplete Market, Kramkov
\cite{Kramkov} constructs the optional decomposition of
supermartingale with respect to a family of equivalent local
martingale measures. He call such a representation optional
because, in contrast to the Doob-Meyer decomposition, it generally
exists only with an adapted (optional) process C. He apply this
decomposition to the problem of hedging European and American
style contingent claims in the setting of incomplete security
markets. Using the optional decomposition, Frey \cite{Frey}
consider construction of superreplication strategies via optimal
stopping which is similar to the optimal stopping problem that
arises in the pricing of American-type derivatives on a family of
probability space with equivalent local martingale measures.

For the realistic financial market, the asset price in the future
is uncertain, the probability distribution of the asset price in
the future is unknown -- which is called Knight uncertain
\cite{Knight}. The probability distribution of the nature state in
the future is unknown, investors have uncertain subjective belief,
which makes their consumption and portfolio choice decisions
uncertain and leads the uncertain asset price in the future.
Pricing contingent claims against such assets under Knight
uncertain is open problem. Peng in \cite{Peng2009} and
\cite{Peng2010} constructs G frame work which is a analysis tool
for nonlinear system and is applied in pricing European contingent
claims under Knight uncertainty \cite{Chen2011}, \cite{Chen2013a}
and \cite{Chen2013b}.

The target of this paper is to establish the bid-ask pricing frame
work for the American contingent claims against risky assets with
G-asset price systems (see \cite{Chen2013b}) on the financial
market under Knight uncertainty. Firstly, on sublinear expectation
space, by using potential theory and sublinear expectation theory
we construct G-Doob-Meyer decomposition for G-supermartingale,
i.e., a right continuous G-supermartingale could be decomposed as
a G-martingale and a right continuous increasing process and the
decomposition is unique. Second, we define bid and ask prices of
the American contingent claim against the assets with G-asset
price systems, and apply the G-Doob-Meyer decomposition to prove
that the bid and ask prices of American contingent claims under
Knight uncertain could be described by the optimal stopping
problems. Finally, we present a free boundary problem, by using
the penalization technique (see Friedman \cite{Friedman}) we
derive that if there exists strong super-solution to the free
boundary problem, then the strong solution to the free boundary
problem exists. And by using truncation and regularization
technique, we prove that the strong solution to the free boundary
problem is the value function of the optimal stopping problem
which is corresponding with pricing problem of the American
contingent claim under Knight uncertain.

The rest of this paper is organized as follows. In Section 2,  we
give preliminaries for the sublinear expectation theory. In
Section 3 we prove G-Doob-Meyer decomposition for
G-supermartingale. In Section 4, using G-Doob-Meyer decomposition,
we construct dynamic superhedge stragies for the optimal stopping
problem, and prove that the solution of the optimal stopping
problem are the bid and ask prices of the American contingent
claims under Knight uncertain. In section 5, we consider a free
boundary problem, prove the strong solution existence of the free
boundary problem, and derive that the solution of the optimal
stopping problem is equivalent the strong solution to the free
boundary problem.

\sce{0}
\section{Preliminaries}

Let $\Omega$ be a given set and let $\cal{H}$ be a linear space of
real valued functions defined on $\Omega$ containing constants.
The space $\cal{H}$ is also called the space of random variables.

\begin{definition}
A sublinear expectation $\hat{E}$ is a functional
$\hat{E}:\mathcal{H}\longrightarrow R$ satisfying

(i) Monotonicity:
$$
\hat{E}[X]\geq \hat{E}[Y]\ \ \mbox{if}\ \ X\geq Y.
$$

(ii) Constant preserving:
$$
\hat{E}[c]=c\ \ \mbox{for}\ \ c\in R.
$$

(iii) Sub-additivity: For each $X,Y\in \cal{H}$,
$$
\hat{E}[X+Y]\leq \hat{E}[X]+\hat{E}[Y].
$$

(iv) Positive homogeneity:
$$
\hat{E}[\lambda X]=\lambda\hat{E}[X]\ \ \mbox{for}\ \ \lambda\geq
0.
$$
The triple $(\Omega,\mathcal{H},\hat{E})$ is called a sublinear
expectation space.
\end{definition}

In this section, we mainly consider the following type of
sublinear expectation spaces $(\Omega,\mathcal{H},\hat{E})$: if
$X_1.X_2,\dots,X_n\in\cal{H}$ then $\varphi(X_1.X_2,\dots,X_n)\in
\cal{H}$ for $\varphi\in C_{b,Lip}(R^n)$, where $C_{b,Lip}(R^n)$
denotes the linear space of functions $\phi$ satisfying
\begin{eqnarray*}
 |\phi(x)-\phi(y)|&\leq& C(1+|x|^m+|y|^m)|x-y| \mbox{ for } x,y\in
 R,\\
&& \mbox{ some } C
> 0, m\in N\mbox{ is depending on }\phi.
\end{eqnarray*}

For each fixed $p\geq 1$, we take $\mathcal{H}_0^p=\{X\in
\mathcal{H},\hat{E}[|X|^p]=0\}$ as our null space, and denote
$\mathcal{H}/\mathcal{H}_0^p$ as the quotient space. We set
$\|X\|_p:=(\hat{E}[|X|^p])^{1/p}$, and extend
$\mathcal{H}/\mathcal{H}_0^p$ to its completion
$\widehat{\cal{H}}_p$ under $\|\cdot\|_p$. Under $\|\cdot\|_p$ the
sublinear expectation $\hat{E}$ can be continuously extended to
the Banach space $(\widehat{\mathcal{H}}_p,\|\cdot\|_p)$. Without
loss generality, we denote the Banach space
$(\widehat{\mathcal{H}}_p,\|\cdot\|_p)$ as
$L^p_G(\Omega,\mathcal{H},\hat{E})$. For the G-frame work, we
refer to \cite{Peng2009} and  \cite{Peng2010}.

In this paper we assume that $\underline{\mu}, \overline{\mu},
\underline{\sigma}$ and $\overline{\sigma}$ are nonnegative
constants such that $\underline{\mu}\leq \overline{\mu}$ and
$\underline{\sigma}\leq\overline{\sigma}$.

\begin{definition} Let $X_1$ and $X_2$ be two random variables in a
sublinear expectation space $(\Omega,\mathcal{H},\hat{E})$, $X_1$
and $X_2$ are called identically distributed, denoted by
$X_1\stackrel{d}{=}X_2$ if
\begin{eqnarray*}
\hat{E}[\phi(X_1)]=\hat{E}[\phi(X_2)]& \mbox{for  } \forall\phi\in
C_{b,Lip}(R^n).
\end{eqnarray*}
\end{definition}
\begin{definition}
In a sublinear expectation space $(\Omega,\mathcal{H},\hat{E})$, a
random variable $Y$ is said to be independent of another random
variable $X$, if
\begin{eqnarray*}
\hat{E}[\phi(X,Y)]=\hat{E}[\hat{E}[\phi(x,Y)]|_{x=X}].
\end{eqnarray*}
\end{definition}
\begin{definition} (G-normal distribution) A random variable $X$
on a sublinear expectation space $(\Omega,\mathcal{H},\hat{E})$ is
called G-normal distributed if
\begin{eqnarray*}
aX+b\bar{X}=\sqrt{a^2+b^2}X&\mbox{for  } a,b\ge 0,
\end{eqnarray*}
where $\bar{X}$ is an independent copy of $X$.
\end{definition}
We denote by $S(d)$ the collection of all $d\times d$ symmetric
matrices. Let $X$ be G-normal distributed random vectors on
$(\Omega,\mathcal{H},\hat{E})$, we define the following sublinear
function
\begin{eqnarray}\label{sublinear-function}
G(A):=\frac{1}{2}\hat{E}[<AX,X>],&A\in S(d).
\end{eqnarray}

\begin{remark}
For a random variable $X$ on the sublinear space
$(\Omega,\mathcal{H},\hat{E})$, there are four typical parameters
to character $X$
\begin{eqnarray*}
\overline{\mu}_X=\hat{E}X,&\underline{\mu}_X=-\hat{E}[-X],\\
\overline{\sigma}_X^2=\hat{E}X^2,&\underline{\sigma}_X^2=-\hat{E}[-X^2],
\end{eqnarray*}
where $[\underline{\mu}_X,\overline{\mu}_X]$ and
$[\underline{\sigma}^2_X,\overline{\sigma}^2_X]$ describe the
uncertainty of the mean and the variance of $X$, respectively.

It is easy to check that if $X$ is G-normal distributed, then
$$
\overline{\mu}_X=\hat{E}X=\underline{\mu}_X=-\hat{E}[-X]=0,
$$
and we denote the G-normal distribution as
$N(\{0\},[\underline{\sigma}^2,\overline{\sigma}^2])$. If $X$ is
maximal distributed, then
$$
\overline{\sigma}_X^2=\hat{E}X^2=\underline{\sigma}_X^2=-\hat{E}[-X^2]=0,
$$
and we denote the maximal distribution (see \cite{Peng2010}) as
$N([\underline{\mu},\overline{\mu}],\{0\})$.
\end{remark}
Let $\mathcal{F}$ as a Borel field subsets of $\Omega$. We are
given a family $\{\mathcal{F}_t\}_{t\in R_+}$ of Borel  subfields
of $\mathcal{F}$, such that
$$
\mathcal{F}_s\subset \mathcal{F}_t,\ \ \ \ s<t.
$$
\begin{definition}
We call $(X_t)_{t\in R}$ a d-dimensional stochastic process on a
sublinear expectation space
$(\Omega,\mathcal{H},\hat{E},\mathcal{F},\{\mathcal{F}\}_{t\in R_+
})$, if for each $t\in R$, $X_t$ is a d-dimensional random vector
in $\cal{H}$.
\end{definition}

\begin{definition}
Let $(X_t)_{t\in R}$ and $(Y_t)_{t\in R}$ be d-dimensional
stochastic processes defined on a sublinear expectation space
$(\Omega,\mathcal{H},\hat{E},\mathcal{F},\{\mathcal{F}\}_{t\in R_+
})$, for each $\underline{t}=(t_1,t_2,\dots,t_n)\in \mathcal{T}$,
$$F_{\underline{t}}^X[\varphi]:=\hat{E}[\varphi(X_{\underline{t}})],\ \ \forall \varphi \in C_{l,Lip}(R^{n\times d})
$$
is called the finite dimensional distribution of $X_t$. $X$ and
$Y$ are said to be identically distributed, i.e., $
X\stackrel{d}{=}Y $, if
$$
F_{\underline{t}}^X[\varphi]=F_{\underline{t}}^Y[\varphi],\ \ \ \
\forall \underline{t}\in \mathcal{T}\ \ \mbox{and}\ \ \forall
\varphi\in C_{l.Lip}(R^{n\times d})
$$
where $\mathcal{T}:=\{\underline{t}=(t_1,t_2,\dots,t_n): \forall
n\in N,t_i\in R,t_i\neq t_j, 0\leq i,j\leq n,i\neq j \}$.
\end{definition}

\begin{definition}\label{Dbm1}
A process $(B_t)_{t\ge 0}$ on the sublinear expectation space
$(\Omega,\mathcal{H},\hat{E},\mathcal{F},\{\mathcal{F}\}_{t\in R_+
})$ is called a G-Brownian motion if the following properties are
satisfied:

(i) $B_0(\omega)=0$;

(ii) For each $t, s>0$, the increment $B_{t+s}-B_{t}$ is G-normal
distributed by
$N(\{0\},[s\underline{\sigma}^2,s\overline{\sigma}^2]$ and is
independent of $(B_{t_1},B_{t_2},\dots,B_{t_n})$, for each $n\in
N$ and $t_1, t_2,\dots,t_n\in (0, t]$;
\end{definition}
From now on, the stochastic process we will consider in the rest
of this paper are all in the sublinear space
$(\Omega,\mathcal{H},\hat{E},\mathcal{F},\{\mathcal{F}\}_{t\in
R_+})$.

\sce{0}
\section{G-Doob-Meyer Decomposition for G-supermartingale}
\begin{definition} A G-supermartingale (resp. G-submartingale) is a real valued process $\{X_t\}$, well adapted to the
$\mathcal{F}_t$ family, such that
\begin{eqnarray}\label{supermartingale}
\begin{array}{lll}
\mbox{(i)}&\ \hat{E}[|X_t|]<\infty & \forall t\in R_+,\\
\mbox{(ii)}&\ \hat{E}[X_{t+s}|\mathcal{F}_s]\leq \mbox{( resp. }
\geq\mbox{ ) } X_s & \forall t\in R_+, \mbox{and }\forall s\in
R_+.
\end{array}
\end{eqnarray}
If equality holds in (ii), the process is a G-martingale.
\end{definition}
We will consider right continuous G-supermartingales, then if
$\{X_t\}$ is right continuous G-supermartingale (ii) in
$(\ref{supermartingale})$ holds with $\mathcal{F}_t$ replaced by
$\mathcal{F}_{t+}$.
\begin{definition}\label{capacity}
Let $A$ be an event in $\mathcal{F}_{t+}$, we define capacity of
$A$ as
\begin{eqnarray}
c(A)=\hat{E}[I_{A}]
\end{eqnarray}
where $I_{A}$ is indicator function of event $A$.
\end{definition}

\begin{definition}
Process $X_t$ and $Y_t$ are adapted to the filtration
$\mathcal{F}_t$. We call $Y_t$ equivalent to $X_t$, if and only if
\begin{eqnarray*}
c(Y_t\neq X_t)=0.
\end{eqnarray*}
\end{definition}

For a right continuous G-supermartingale $\{X_t\}$ with
$\hat{E}[X_t]$ is right continuous function of t, we can find a
right continuous G-supermartingale $\{Y_t\}$ equivalent to
$\{X_t\}$ by define
\begin{eqnarray*}
Y_t(\omega):=X_{t+}(\omega)=\lim_{s\downarrow t}X_s(\omega),
&\mbox{for any } \omega\in \Omega\ \ \mbox{such that the limit
exits}\\
Y_t(\omega):=0,&\mbox{otherwise}.
\end{eqnarray*}
Without loss generality, we denote
$\mathcal{F}_t=\mathcal{F}_{t+}.$
\begin{definition}
For a positive constant $T$, we define stop time $\tau$ in $[0,T]$
as a positive, random variable $\tau(\omega)$ such that,
$\{\tau\leq T\}\in \mathcal{F}_T$.
\end{definition}
Let $\{X_t\}$ be a right continuous G-supermartingale, denote
$X_{\infty}$ as the last element of the process $X_t$, then the
process $\{X_t\}_{0\leq t\leq \infty}$ is a G-supermartingale.

\begin{definition}
A right continuous increasing process is a well adapted stochastic
process $\{A_t\}$ such that:

(i) $A_0=0$ a.s.

{ii} For almost every $\omega$, the function $t\longrightarrow
A_t(\omega)$ is positive, increasing, and right continuous. Let
$A_{\infty}(\omega):=\lim_{t\longrightarrow \infty}A_t(\omega)$,
we shall say that the right continuous increasing process is
integrable if $\hat{E}[A_{\infty}]<\infty$.
\end{definition}
\begin{definition}\label{de-natural}
An increasing process $A$ is called natural if for every bounded,
right continuous G-martingale $\{M_t\}_{0\leq t< \infty }$ we have
\begin{eqnarray}\label{eq-natural}
\hat{E}[\int_{(0,t]}M_sdA_s]=\hat{E}[\int_{(0,t]}M_{s-}dA_s],&\mbox{for
every  } 0<t<\infty
\end{eqnarray}
\begin{lemma}\label{le-natural} If $A$ is an increasing process and $\{M_t\}_{0\leq
t<\infty}$ is bounded, right continuous G-martingale, then
\begin{eqnarray}\label{la-eq-natural1}
\hat{E}[M_tA_t]=\hat{E}[\int_{(0,t]}M_{s}dA_s].
\end{eqnarray}
In particular, condition $(\ref{eq-natural})$ in Definition
$\ref{de-natural}$ is equivalent to
\begin{eqnarray}\label{la-eq-natural2}
\hat{E}[M_tA_t]=\hat{E}[\int_{(0,t]}M_{s-}dA_s].
\end{eqnarray}
\end{lemma}
{\bf Proof. } For a partition $\Pi=\{t_0,t_1,\cdots,t_n\}$ of
$[0,t]$, with $0=t_0\leq t_1\leq \cdots\leq t_n=t$, we define
$$
M_s^{\Pi}=\sum_{k=1}^{n}M_{t_n}I_{(t_{k-1},t_k]}(s).
$$
Since $M$ is G-martingale
\begin{eqnarray*}
\hat{E}[\int_{(0,t]}M_{s}^{\Pi}dA_s]&=&\hat{E}[\sum_{k=1}^{n}M_{t_{k}}(A_{t_k}-A_{t_{k-1}})]\\
&=&\hat{E}[\sum_{k=1}^{n}M_{t_{k}}A_{t_k}-\sum_{k=1}^{n-1}M_{t_{k+1}}A_{t_k}]\\
&=&\hat{E}[M_{t}A_{t}-\sum_{k=1}^{n-1}(M_{t_{k+1}}-M_{t_k})A_{t_k}]\\
&=&\hat{E}[M_{t}A_{t}-\sum_{k=1}^{n-1}(M_{t_{k+1}}-M_{t_k})A_{t_k}]\\
&=&\hat{E}[M_{t}A_{t}],
\end{eqnarray*}
we finish the proof of the Lemma.$\ \ \ \ \square$

\end{definition}
\begin{definition} A positive right continuous G-supermartingale
$\{Y_t\}$ with $\lim_{t\longrightarrow \infty}Y_t(\omega)=0$ is
called a potential.
\end{definition}

\begin{definition}
For $a\in [0,\infty]$, a process $\{X_t, t\in [0,a]\}$is said to
be uniformly integrable on $[0,a]$ if
$$
\sup_{t\in[0,a]}\hat{E}[|X_t|I_{|X_t|>x}]\longrightarrow 0,\ \
\mbox{as } x\longrightarrow 0.
$$
\end{definition}
\begin{definition}
Let $a\in [0,\infty]$, and let $\{X_t\}$ be a right continuous
process,
we shall say that it belongs to the class (GD) on this interval,
if all the random variable $X_T$ are uniformly integrable, $T
\mbox{be stop time bounded by }a$. If $\{X_t\}$ belongs to the
class (GD) on every interval $[0,a],a<\infty$, it will be said to
belong locally to the class (GD).
\end{definition}

If $\{A_t\}$ is a integrable right continuous, increasing process,
then the process $\{-A_t\}$ is a negative G-supermartingale, and
$\{\hat{E}[A_{\infty}|\mathcal{F}_t]-A_t\}$ is a potential of the
class (GD), which we shall call the potential generated by
$\{A_t\}$.

\begin{proposition}

(1) Any right continuous G-martingale $\{X_t\}$ belongs locally to
the class (GD).

(2) Any right continuous G-supermartingale $\{X_t\}$, which is
bounded from above, belongs locally to the class (GD).

(3) Any right continuous supermartingale $\{X_t\}$, which belongs
locally to the class (GD) and is uniformly integrable, belongs to
the class (GD).
\end{proposition}
{\bf Proof.} (1) If $a<\infty$, and $T$ is a stop time, $T\leq a$,
then G-martingale process $\{X_t\}$ has
$X_T=\hat{E}[X_a|\mathcal{F}_T]$. Hence
\begin{eqnarray*}
\hat{E}[|X_T|I_{\{|X_T|>n\}}]\leq \hat{E}[|X_a|I_{\{|X_T|>n\}}]
\end{eqnarray*}
As $n\cdot c(|X_T|>n)\leq \hat{E}[|X_T|]\leq \hat{E}[|X_{a}|]$, we
have $c(|X_T|>n)\longrightarrow 0$ as $n\longrightarrow \infty$,
then $\hat{E}[|X_a|I_{\{|X_T|>n\}}]\leq
(\hat{E}[|X_d|^2])^{1/2}(c(I_{\{|X_T|>n\}}))^{1/2} \longrightarrow
0$ as $n\longrightarrow \infty$, from which we prove (1).

(2) If $a<\infty$, and $T$ is a stop time, $T\leq a$, then
G-supermartingale process $\{X_t\}$ has $X_T\geq
\hat{E}[X_a|\mathcal{F}_T]$. Suppose that $\{X_t\}$ is negative,
then
\begin{eqnarray*}
\hat{E}[-X_TI_{\{X_T<-n\}}]\leq \hat{E}[-X_aI_{\{X_T<-n\}}]
\end{eqnarray*}
we complete the proof of (2) by using the similar argument in
proof (1).

(3) $\{X_t\}$ is uniformly integrable, we set
\begin{eqnarray*}
X_t=\hat{E}[X_{\infty}|\mathcal{F}_t]+(X_t-\hat{E}[X_{\infty}|\mathcal{F}_t]).
\end{eqnarray*}

The fist part on the right hand of the above equation
$\hat{E}[X_{\infty}|\mathcal{F}_t]$ is a G-martingale, and
equivalent to a right continuous process, and from (1) we know
that it belongs to the class (GD). We denote the second part in
the above equation as $\{Y_t\}$, it is a potential, i.e., a
positive right continuous G-supermartingale, and
$\lim_{t\longrightarrow \infty}Y_t(\omega)=0$ a.s.. Next we will
prove that $\{Y_t\}$ belongs to the class (GD). Since both
$\inf(T,a)$ and $\sup(T,a)$ are stop times
\begin{eqnarray*}
\hat{E}[Y_TI_{\{Y_T>n\}}]&\leq& \hat{E}[Y_TI_{\{T\leq a,
Y_T>n\}}]+\hat{E}[Y_TI_{\{T>a\}}]\\
&\leq& \hat{E}[Y_TI_{\{T\leq a, Y_T>n\}}]+\hat{E}[Y_a].
\end{eqnarray*}
As $\lim_{a\longrightarrow \infty}\hat{E}[Y_a]=0$ and $\{Y_t\}$
locally belongs to (GD), i.e., $\lim_{n\longrightarrow \infty
}\hat{E}[Y_TI_{\{T\leq a, Y_T>n\}}]=0$, which prove that
\begin{eqnarray*}
\lim_{n\longrightarrow \infty}\hat{E}[Y_TI_{\{Y_T>n\}}]=0.
\end{eqnarray*}
We complete the proof.$\ \ \square$

\begin{lemma} \label{XMA}Let $\{X_t\}$ be a right continuous
G-supermartingale, and $\{X_t^n\}$ a sequence of decomposed right
continuous G-supermartingale:
$$
X_t^n=M_t^n-A_t^n,
$$
where $\{M_t^n\}$ is G-martingale, and $\{A_t^n\}$ is right
continuous increasing process. Suppose that, for each $t$, the
$X_t^n$ converge to $X_t$ in the $L_G^1(\Omega)$ topology, and the
$A_t^n$ are uniformly integrable in $n$. Then the decomposition
problem is solvable for the G-supermartingale $\{X_t\}$, more
precisely, there are a right continuous increasing process
$\{A_t\}$, and a G-martingale $\{M_t\}$, such that $X_t=M_t-A_t$.
\end{lemma}
{\bf Proof. } We denote by $w$ the weak topology
$w(L^1_G(\Omega),L^{\infty}_G(\Omega))$, a sequence of integrable
random variables $f_n$ converges to a random variable $f$ in the
$w$-topology, if and only if $f$ is integrable, and
\begin{eqnarray*}
\lim_{n\longrightarrow \infty}\hat{E}[f_ng]=\hat{E}[fg],&\forall
g\in L^{\infty}_G(\Omega).
\end{eqnarray*}
Since the $A_t^n$ are uniformly integrable in $n$, by the
properties of the sublinear expectation $\hat{E}[\cdot]$ there
exists a $w$-convergent subsequence $A_t^{n_k}$ converge in the
$w$-topology to the random variables $A_t^{\prime}$, for all
rational values of $t$. To simplify the notations, we shall use
$A_t^{n}$ converge to $A_t^{\prime}$ in the $w$-topology for all
rational values of $t$. An integrable random variable $f$ is
$\mathcal{F}_t$-measurable if and only if it is orthogonal to all
bounded random variables $g$ such that
$\hat{E}[g|\mathcal{F}_t]$=0, it follows that $A_t^{\prime}$ is
$\mathcal{F}_t-$measurable. For $s<t$, $s$ and $t$ rational,
\begin{eqnarray}
\hat{E}[(A_t^{n}-A_s^{n})I_{B}]\geq 0
\end{eqnarray}
where $B$ denote any $\mathcal{F}$ set.

As $X_t^n$ converge to $X_t$ in $L^{1}_G(\Omega)$ topology, which
is in a stronger topology than $w$, the $M_t^n$ converge to random
variables $M_t^{\prime}$ for $t$ rational, and the process
$\{M_t^{\prime}\}$ is G-martingale; then there is a right
continuous G-martingale $\{M_t\}$, defined for all values of $t$,
such that $c(M_t\neq M_t^{\prime})=0$ for each rational $t$. We
define $A_t=X_t+M_t$, $\{A_t\}$ is a right continuous increasing
process, or at least becomes so after a modification on a set of
measure zero. We complete the proof. $\ \ \square$
\begin{lemma} \label{potentialAH}
Let $\{X_t\}$ be a potential and belong to the class (GD). We
consider the measurable, positive and well adapted processes
$H=\{H_t\}$ with the property that the right continuous increasing
processes
$$
A(H)=\{A_t(H,\omega)\}=\{\int_0^tH_s(\omega)ds\}
$$
are integrable, and the potentials $Y(H)=\{Y_t(H,\omega)\}$ they
generate are majorized by $X_t$. Then, for each $t$, the random
variables $A_t(H)$ of all such processes $A(H)$ are uniformly
integrable.
\end{lemma}
{\bf Proof. } It is sufficient to prove that the $A_{\infty}(H)$
are uniformly integrable.

(1) First we assume that $X_t$ is bounded by some positive
constant $C$, then $\hat{E}[A_{\infty}^{2}(H)]\leq 2 C^2$, and the
uniform integrability follows.

We have that
\begin{eqnarray*}
A_{\infty}^2(H,\omega)&=&2\int_0^{\infty}[A_{\infty}(H,\omega)-A_{u}(H,\omega)]dA_u(H,\omega)\\
&=&2\int_0^{\infty}[A_{\infty}(H,\omega)-A_{u}(H,\omega)]H_u(\omega)du.
\end{eqnarray*}
By using the sub-additive property of the sublinear expectation
$\hat{E}$
\begin{eqnarray*}
\hat{E}[A_{\infty}^2(H,\omega)]&=&\hat{E}[\hat{E}A_{\infty}^2(H,\omega)|\mathcal{F}_t]]\\
&\leq &2\hat{E}[\int_0^{\infty}H_u\hat{E}[A_{\infty}(H,\omega)-A_{u}(H,\omega)|\mathcal{F}_u]du]\\
&=&2\hat{E}[\int_0^{\infty}H_uY_u(H)du]\\
& \leq & 2C\hat{E}[\int_0^{\infty}H_udu]\\
&=&2C\hat{E}[Y_0(H)]\\
&\leq & 2C^2.
\end{eqnarray*}

(2) In order to prove the general case, it will be enough to prove
that any $H$ such that $Y(H)$ is majorized by $\{X_t\}$ is equal
to a sum $H^c+H_c$, where (i) $A(H^c)$ generates a potential
bounded by $c$, and (ii) $\hat{E}[A_{\infty}[H_c]]$ is smaller
than some number $\varepsilon_c$, independent of $H$, such that
$\varepsilon_c\longrightarrow 0$ as $c\longrightarrow 0$. Define
$$
H_t^c(\omega)=H_t(\omega)I_{\{X_t(\omega)\in [0,c]\}},\ \
H_{ct}=H_t-H_t^c.
$$
Set
$$
T^c(\omega)=\inf\{t:\ \mbox{such that }X_t(\omega)\geq c\},
$$
as $c$ goes to infinity $\lim_{c\longrightarrow
\infty}T^c(\omega)=\infty$, therefore $X_{T^c}\longrightarrow 0$,
and the class (GD) property implies that
$\hat{E}[X_{T^c}]\longrightarrow 0$. $T^c$ is a stop time, and
$I_{\{X_t(\omega)\in [0,c]\}}=1$ before time $T^c$. Hence
\begin{eqnarray*}
\hat{E}[A_{\infty}(H_c)]&=&\hat{E}[\int_0^{\infty}H_u(1-I_{\{X_u(\omega)\in
[0,c]\}})]du\\
&\leq &\hat{E}[\int_0^{\infty}H_udu]\\
&=&\hat{E}[A_{\infty}(H)-A_{T^c}(H)]\\
&=&\hat{E}[\hat{E}[A_{\infty}(H)-A_{T^c}(H)|\mathcal{F}_t]]\\
&=&\hat{E}[Y_{T^c}(H)]\leq \hat{E}[X_{T^c}(H)]\\
&\leq &\varepsilon_c,\ \ \mbox{for large enough\ } c,
\end{eqnarray*}
from which we prove (ii). We shall prove (i), first we prove that
$Y(H^c)$ is bounded by $c$:
\begin{eqnarray}\label{Yt}
\begin{array}{rcl}
Y_t(H^c)&=&\hat{E}[A_{\infty}(H^c)-A_t(H^c)|\mathcal{F}_t]\\
&=&\hat{E}[\int_t^{\infty}H_uI_{\{X_u(\omega)\in
[0,c]\}}du|\mathcal{F}_t]\\
&\leq &\hat{E}[\int_{S^c}^{\infty}H_uI_{\{X_u(\omega)\in
[0,c]\}}du|\mathcal{F}_t]\\
&=&\hat{E}[\hat{E}[\int_{S^c}^{\infty}H_uI_{\{X_u(\omega)\in
[0,c]\}}du|\mathcal{F}_{S^c}]|\mathcal{F}_t]\\
&=&\hat{E}[Y_{S^c}|\mathcal{F}_{t}]\\
&\leq & c,
\end{array}
\end{eqnarray}
where we set
$$
S^c(\omega)=\inf\{t:\ \mbox{such that }X_t(\omega)\leq c\},
$$
and use
$$
\int_t^{S^c(\omega)}H_uI_{\{X_u(\omega)\in [0,c]\}}du=0.
$$
the inequality $(\ref{Yt})$ holds for each $t$, therefore for
every rational $t$, and for every $t$ in consideration of the
right continuity, which complete the proof. $\ \ \square$

\begin{lemma}\label{pXt} Let $\{X_t\}$ be a potential and belong to the class (GD), $k$ is a positive number,
define $Y_t=\hat{E}[X_{t+k}|\mathcal{F}_t]$, then $\{Y_t\}$ is a
G-supermartingale. Denote by $\{p_kX_t\}$ a right continuous
version of $\{Y_t\}$, then $\{p_kX_t\}$ is potential.

Use the same notations as in Lemma $\ref{potentialAH}$. Let $k$ be
a positive number, and
$H_{t,k}(\omega)=(X_t(\omega)-p_kX_t(\omega))/k$. The process
$H_k=\{H_{t,k}\}$ verify the assumptions of Lemma
$\ref{potentialAH}$, and their potentials increase to $\{X_t\}$ as
$k\longrightarrow 0$.
\end{lemma}
{\bf Proof. } If $t<u$
\begin{eqnarray*}
&&\hat{E}[\frac{1}{k}(\int_0^u[X_s-p_kX_s]ds-\int_0^t[X_s-p_kX_s]ds)|\mathcal{F}_t]\\
&=&\hat{E}[\frac{1}{k}\int_t^u[X_s-p_kX_s]ds|\mathcal{F}_t].
\end{eqnarray*}
For $s\geq t$, $\hat{E}[p_k X_s|\mathcal{F}_t]=\hat{E}[\hat{E}
[X_{s+k}|\mathcal{F}_s]|\mathcal{F}_t]=\hat{E}[X_{s+k}|\mathcal{F}_t].$
We have that
\begin{eqnarray*}
\hat{E}[\frac{1}{k}\int_t^u[X_s-p_kX_s]ds|\mathcal{F}_t]\geq
\hat{E}[\frac{1}{k}\int_t^{t+k}X_sds|\mathcal{F}_t]-\hat{E}[\frac{1}{k}\int_u^{u+k}X_sds|\mathcal{F}_t],
\end{eqnarray*}
by the sub-additive property of the sublinear expectation
$\hat{E}$, we derive that
\begin{eqnarray*}
&&\hat{E}[\frac{1}{k}\int_t^{t+k}X_sds|\mathcal{F}_t]-\hat{E}[\frac{1}{k}\int_u^{u+k}X_sds|\mathcal{F}_t]\\
&\geq
&\hat{E}[\frac{1}{k}\int_t^{t+k}X_sds|\mathcal{F}_t]-\frac{1}{k}\int_u^{u+k}\hat{E}[X_s|\mathcal{F}_t]ds\\
&\geq & \hat{E}[\frac{1}{k}\int_t^{t+k}X_sds|\mathcal{F}_t]-X_t\\
&\geq & -\hat{E}[\frac{1}{k}\int_t^{t+k}(X_t-X_s)ds|\mathcal{F}_t]\\
&\geq & -\frac{1}{k}\int_t^{t+k}\hat{E}[X_t-X_s|\mathcal{F}_t]ds\\
&\geq & 0.
\end{eqnarray*}
Hence, we derive that, for any $u,t$ such that $u>t$
\begin{eqnarray*}
\hat{E}[\frac{1}{k}\int_t^u[X_s-p_kX_s]ds|\mathcal{F}_t]\geq 0.
\end{eqnarray*}
If there exits $s_0\geq 0$ such that
$\frac{1}{k}[X_{s_0}-p_kX_{s_0}]<0$, the right continuous of
$\{X_t\}$ implies that there exists $\delta>0$ such that
$\frac{1}{k}[X_{s}-p_kX_{s}]<0$ on the interval
$[s_0,s_0+\delta]$. Thus
\begin{eqnarray*}
\hat{E}[\frac{1}{k}\int_{s_0}^{s_0+\delta}[X_s-p_kX_s]ds|\mathcal{F}_{s_0}]<0,
\end{eqnarray*}
which is contradiction, we prove that
$(X_t(\omega)-p_kX_t(\omega))/k$ is a positive, measurable and
well adapted process.

Since $\{X_t\}$ is right continuous G-supermartingale
$$\lim_{s\downarrow t}X_s=X_t $$
\begin{eqnarray*}
\lim_{k\downarrow 0}Y_t(H_k)&=&\lim_{k\downarrow 0}\hat{E}[\frac{1}{k}\int_t^{\infty}[X_s-p_kX_s]ds|\mathcal{F}_t]\\
&=&\lim_{k\downarrow 0}\hat{E}[\frac{1}{k}\int_t^{t+k}X_sds|\mathcal{F}_t]\\
&=&\hat{E}[\lim_{k\downarrow 0}\frac{1}{k}\int_t^{t+k}X_sds|\mathcal{F}_t]\\
&=&X_t,
\end{eqnarray*}
we finish the proof.$\ \ \ \ \square$

From Lemma $\ref{XMA}$, $\ref{potentialAH}$, and $\ref{pXt}$ we
can prove the following Theorem
\begin{theorem}\label{potentialdecomposition}
A potential $\{X_t\}$ belongs to the class (GD) if, and only if,
it is generated by some integrable right continuous increasing
process.
\end{theorem}

\begin{theorem}\label{GDoob} (G-Doob-Meyer's Decomposition)

(1) $\{X_t\}$ is a right continuous G-supermartingale if and only
if it belongs to the class (GD) on every finite interval. More
precisely, $\{X_t\}$ is then equal to the difference of a
G-martingal $M_t$ and a right continuous increasing process $A_t$
\begin{eqnarray}
X_t=M_t-A_t.
\end{eqnarray}

(2) If the right continuous increasing process $A$ is natural, the
decomposition is unique.
\end{theorem}
{\bf Proof. } (1) The necessity is obvious. We will prove the
sufficiency, we choose a positive number $a$ and define
\begin{eqnarray}
X_t^{\prime}(\omega):=X_t(\omega),\ \ t\in [0,a]&
X_t^{\prime}(\omega):=X_a(\omega),\ \ t>a,
\end{eqnarray}
the $\{X_t^{\prime}\}$ is a right continuous G-supermartingale of
the class (GD), by Theorem $\ref{potentialdecomposition}$ there
exist the following decomposition
$$
X_t^{\prime}=M_t^{\prime}-A_t^{\prime},
$$
where $\{M_t^{\prime}\}$ is a G-martingal, and $\{A_t^{\prime}\}$
is a right continuous increasing process.

Let $a\longrightarrow \infty$, as in Lemma $\ref{pXt}$ the
expression of the $Y_t(H_k)$ that the $A_t^{\prime}$ depend only
on the values of $\{X_t^{\prime}\}$ on intervals
$[0,t+\varepsilon]$, with $\varepsilon$ small enough. As
$a\longrightarrow \infty$, they don't vary any more once $a$ has
reached values greater than $t$, us again Lemma $\ref{XMA}$, we
finish the proof of the Theorem.

(2) Assume that $X$ admits both decompositions
$$
X_t=M_t^{\prime}-A_t^{\prime}=M_t^{\prime\prime}-A_t^{\prime\prime},
$$
where $M_t^{\prime}$ and $M_t^{\prime\prime}$ are G-martingale and
$A_t^{\prime}$, $A_t^{\prime\prime}$ are natural increasing
process. We define
$$
\{C_t:=A_t^{\prime}-A_t^{\prime\prime}=M_t^{\prime}-M_t^{\prime\prime}\}.
$$
Then $\{C_t\}$ is a G-martingale, and for every bounded and right
continuous G-martingale $\{\xi_t\}$, from Lemma $\ref{le-natural}$
we have
$$
\hat{E}[\xi_t(A_t^{\prime}-A_t^{\prime\prime})]=\hat{E}[\int_{(0,t]}\xi_{s-}dC_s=\lim_{n\longrightarrow
\infty}\sum_{k=1}^{m_n}\xi_{t_{j-1}^{n}}[C_{t_{j}^{(n)}}-C_{t_{j-1}^{(n)}}],
$$
where $\Pi_n=\{t_0^{(n)},\cdots,t_{m_n}^{(n)}\},n\geq 1$ is a
sequence of partitions of $[0,t]$ with $\max_{1\leq j\leq
m_n}(t_{j}^{(n)}-t_{j-1}^{(n)})$ converging to zero as
$n\longrightarrow \infty$. Since $\xi$ and $C$ are both
G-martingale, we have
$$
\hat{E}[\xi_{t_{j-1}^{(n)}}(C_{t_{j}^{(n)}}-C_{t_{j-1}^{(n)}})]=0,\
\ \mbox{and thus }\ \
\hat{E}[\xi_{t_{j-1}}(A_{t}^{\prime}-A_{t}^{\prime\prime})]=0.
$$
For an arbitrary bonded random variable $\xi$, we can select
$\{\xi_t\}$ to be a right-continuous equivalent process of
$\{\hat{E}[\xi|\mathcal{F}_t]\}$, we obtain that
$\hat{E}[\xi(A_{t}^{\prime}-A_{t}^{\prime\prime})]=0$. We set
$\xi=I_{A_{t}^{\prime}\neq A_{t}^{\prime\prime}}$ therefore
$c(A_{t}^{\prime}\neq A_{t}^{\prime\prime})=0.$
 $\ \ \ \ \square$

By Theorem $\ref{GDoob}$ and G-martingale decomposition Theorem in
\cite{Peng2010} and \cite{Song}, we have the following
G-Doob-Meyer's Theorem
\begin{theorem}\label{GDoobB}
$\{X_t\}$ is a right continuous G-supermartingale, there exists a
right continuous increasing process $A_t$ and adapted process
$\eta_t$, such that
\begin{eqnarray}
X_t=\int_0^t\eta_s dB_s-A_t,
\end{eqnarray}
where $B_t$ is G-Brownian motion.
\end{theorem}

\sce{0}
\section{Superhedging strategies and optimal stopping}
\subsection{Financial model and G-asset price system}
We consider a financial market with a nonrisky asset (bond) and a
risky asset (stock) continuously trading in market. The price P(t)
of the bond is given by
\begin{eqnarray}\label{bond}
dP(t) = rP(t)dt, P(0) = 1,
\end{eqnarray}
where $r$ is the short interest rate, we assume a constant
nonnegative short interest rate. We assume the risk asset with the
G-asset price system $((S_{u})_{u\geq t},\hat{E})$ (see
\cite{Chen2013b}) on sublinear expectation space
$(\Omega,\mathcal{H},\hat{E},\mathcal{F},(\mathcal{F}_t))$ under
Knightian uncertainty, for given $t\in [0,T]$ and $x\in R^d$
\begin{eqnarray}\label{G-asset}
\begin{array}{l}
dS_{u}^{t,x}= S_{u}^{t,x}dB_t=S_{u}^{t,x}(db+d\hat{B}_t),\\
S_{t}^{t,x}=x
\end{array}
\end{eqnarray}
where $B_t$ is the generalized G-Brownian motion. The uncertain
volatility is described by the G-Brownian motion $\hat{B}_t$. The
uncertain drift $b_t$ can be rewritten as
$$
b_t=\int_0^t\mu_tdt
$$
where $\mu_t$ is the asset return rate (\cite{Chen2011}). Then the
uncertain risk premium of the G-asset price system
\begin{eqnarray}
\theta_t=\mu_t-r,
\end{eqnarray}
is uncertain and distributed by
$N([\underline{\mu}-r,\overline{\mu}-r],\{0\})$ (\cite{Chen2011}),
where $r$ is the interest rate of the bond.

Define
\begin{eqnarray}\label{GGirsanov}
\widetilde{B}_t:=B_t-rt=b_t+\hat{B}_t-rt,
\end{eqnarray}
we have the following G-Girsanov Theorem (presented in
\cite{Chen2013a}, \cite{Chen2013b} and \cite{Humingshang2012})
\begin{theorem} {\bf (G-Girsanov Theorem)} Assume that $(B_t)_{t\geq 0}$
is generalized G-Brownian motion on
$(\Omega,\mathcal{H},\hat{E},\mathcal{F}_t)$, and
$\widetilde{B}_t$ is defined by $(\ref{GGirsanov})$, there exists
G-expectation space $(\Omega,\mathcal{H},E^G,\mathcal{F}_t)$ such
that $\widetilde{B}_t$ is G-Brownian motion under the G-
expectation $E^G$, and
\begin{eqnarray}
\hat{E}[\hat{B}_t^2]=E^G[\tilde{B}_t^2],&-\hat{E}[-\hat{B}_t^2]=-E^G[-\tilde{B}_t^2].
\end{eqnarray}
\end{theorem}
By the G-Girsanov Theorem, the G-asset price system
($\ref{G-asset}$) of the risky asset can be rewritten on
$(\Omega,\mathcal{H},E^G,\mathcal{F}_t)$ as follows
\begin{eqnarray}\label{G-assetneutral}
\begin{array}{l}
dS_{u}^{t,x}=S_{u}^{t,x}(rdt+d\tilde{B}_t),\\
S_{t}^{t,x}=x,
\end{array}
\end{eqnarray}
then by G-Ito formula we have
\begin{eqnarray}
S_{u}^{t,x}=x\exp{(r(u-t)+\tilde{B}_{u-t}-\frac{1}{2}(<\tilde{B}_u>-<\tilde{B}_t>))},u>t
\end{eqnarray}

\subsection{Construction of superreplication strategies  via
optimal stopping}

We consider the following class of contingent claims:
\begin{definition}\label{payoff}We define a class of contingent claims with the nonnegative payoff $\xi\in L^2_G (\Omega_T )$
 has the following form
\begin{eqnarray}
\xi=f(S_T^{t,x})
\end{eqnarray}
for some function $f:\Omega\longrightarrow R$ such that the
process
\begin{eqnarray}
f_u:= f(S_u^{t,x})
\end{eqnarray}
is bounded below and c$\grave{a}$dl$\grave{a}$g.
\end{definition}
We consider a contingent claim $\xi$ with payoff defined in
Definition $\ref{payoff}$ written on the stockes $S_t$ with
maturity $T$. We give definitions of superhedging (resp.
subhedging) strategy and ask (resp. bid) price of the claim $\xi$.
\begin{definition} (1) A self-financing superstrategy (resp.
substrategy) is a vector process $(Y,\pi,C)$ (resp. $(-Y,\pi,C)$),
where $Y$ is the wealth process, $\pi$ is the portfolio process,
and $C$ is the cumulative consumption process, such that
\begin{eqnarray}
 dY_t =rY_tdt +\pi_td\tilde{B}_t-dC_t,\\
 \mbox{(resp. } -dY_t = -rY_tdt +\pi_td\tilde{B}_t-dC_t\mbox{ )}
 \end{eqnarray}
 where C is an increasing, right-continuous
process with $C_0 = 0$. The superstrategy (resp. substrategy) is
called feasible if the constraint of nonnegative wealth holds
$$
Y_t\geq 0,\ \ t\in[0,T].
$$

(2) A superhedging (resp. subhedging) strategy against the
European contingent claim $\xi$ is a feasible self-financing
superstrategy $(Y,\pi,C)$ (resp. substrategy $(-Y,\pi,C)$) such
that $Y_T = \xi$ (resp. $-Y_T =-\xi$). We denote by $\mathcal{H}
(\xi)$ (resp. $\mathcal{H}^{\prime}(-\xi)$) the class of
superhedging (resp. subhedging) strategies against $\xi$, and if
$\mathcal{H} (\xi)$ (resp. $\mathcal{H}^{\prime}(-\xi)$) is
nonempty, $\xi$ is called superhedgeable (resp. subhedgeable).

(3) The ask-price $X(t)$ at time $t$ of the superhedgeable claim
$\xi$ is defined as
$$
X(t)=\inf\{x\ge 0:\exists (Y_t,\pi_t,C_t)\in\mathcal{H}(\xi)\mbox{
such that } Y_t=x\},
$$
and bid-price $X^{\prime}(t)$ at time $t$ of the subhedgeable
claim $\xi$ is defined as
$$
X^{\prime}(t)=\sup\{x\ge 0:\exists
(-Y_t,\pi_t,C_t)\in\mathcal{H}^{\prime}(-\xi)\mbox{ such that }
-Y_t=-x\}.
$$
\end{definition}

Under uncertainty, the market is incomplete and the superhedging
(resp. subhedging) strategy of the claim is not unique. The
definition of the ask-price $X(t)$ implies that the ask-price
$X(t)$ is the minimum amount of risk for the buyer to superhedging
the claim, then it is coherent measure of risk of all
superstrategies against the claim for the buyer. The coherent risk
measure of all superstrategies against the claim can be regard as
the sublinear expectation of the claim, we have the following
representation of bid-ask price of the claim via optimal stopping
(Theorem $\ref{th-bid-ask}$) .

Let $(\mathcal{G}_t)$ be a filtration on G-expectation space
$(\Omega,\mathcal{H},E^G,\mathcal{F},(\mathcal{F}_t)_{t\geq 0})$,
and $\tau_1$ and $\tau_2$ be $(\mathcal{G}_t)$- stopping times
such that $\tau_1\leq \tau_2$ a.s.. We denote by
$\mathcal{G}_{\tau_1,\tau_2}$ the set of all finite
$(\mathcal{G}_t)$-stopping times $\tau$ with $\tau_1\leq \tau\leq
\tau_2$.

For given $t\in [0,T]$ and $x\in R_+^d$, we define the function
$V^{Am}:[0,T]\times\Omega\longrightarrow R$ as the value function
of the following optimal-stopping problem
\begin{eqnarray}\label{valuefunction}
V^{Am}(t,S_t)&:=&\sup_{\nu\in
\mathcal{F}_{t,T}}E^G_t[f_{\nu}]\\
&=&\sup_{\nu\in \mathcal{F}_{t,T}}E^G_t[f(S_{\nu})]
\end{eqnarray}

\begin{proposition}\label{p-optimaltime}Consider two stopping times $\underline{\tau}\leq
\overline{\tau}$ on filtration $\mathcal{F}$. Let $(f_t)_{t\geq
0}$ denote some adapted and RCLL-stochastic process, which is
bounded below. Then we have for two points $s, t\in
[0,\overline{\tau}]$ and $s<t$
\begin{eqnarray}
\mbox{ess
sup}_{\tau\in\mathcal{F}_{\underline{\tau},\overline{\tau}}}\{E^G_s[f_{\tau}]\}
=E^G_s[\mbox{ess
sup}_{\tau\in\mathcal{F}_{\underline{\tau},\overline{\tau}}}\{E^G_t[f_{\tau}]\}]
\end{eqnarray}
\end{proposition}
{\bf Proof.} By the consistent property of the conditional
G-expectation, for $\tau\in
\mathcal{F}_{\underline{\tau},\overline{\tau}}$, $s, t\in
[0,\overline{\tau}]$ and $s<t$
\begin{eqnarray*}
E^G_s[f_{\tau}]&=&E^G_s[E^G_t[f_{\tau}]]\\
&\leq & E^G_s[\mbox{ess
sup}_{\tau\in\mathcal{F}_{\underline{\tau},\overline{\tau}}}\{E^G_t[f_{\tau}]\}],
\end{eqnarray*}
thus we have
\begin{eqnarray*}
\mbox{ess
sup}_{\tau\in\mathcal{F}_{\underline{\tau},\overline{\tau}}}\{E^G_s[f_{\tau}]\}
\leq E^G_s[\mbox{ess
sup}_{\tau\in\mathcal{F}_{\underline{\tau},\overline{\tau}}}\{E^G_t[f_{\tau}]\}].
\end{eqnarray*}
There exists a sequences $\{\tau_n\}\longrightarrow \tau^*\in
[\underline{\tau},\overline{\tau}]$ as $n\longrightarrow \infty$,
such that
\begin{eqnarray}
\lim_{n\longrightarrow \infty}E^G_t[f_{\tau_n}]=
E^G_t[f_{\tau^*}]= \mbox{ess
sup}_{\tau\in\mathcal{F}_{\underline{\tau},\overline{\tau}}}\{E^G_t[f_{\tau}]\},
\end{eqnarray}
notice that
\begin{eqnarray*}
&&E^G_s[\mbox{ess
sup}_{\tau\in\mathcal{F}_{\underline{\tau},\overline{\tau}}}\{E^G_t[f_{\tau}]\}]\\
&=&E_s^G[E^G_t[f_{\tau^*}]]\\
&=&E_s^G[f_{\tau^*}]\\
&\leq & \mbox{ess
sup}_{\tau\in\mathcal{F}_{\underline{\tau},\overline{\tau}}}\{E^G_s[f_{\tau}]\},
\end{eqnarray*}
we prove the Proposition. $\ \ \ \ \ \ \square$
\begin{proposition} \label{p-supermartingale}
The process $V^{Am}(t,S_t)_{0\leq t\leq T}$ is a G-supermartingale
in $(\Omega,\mathcal{H},E^G,\mathcal{F},\mathcal{F}_t)$.
\end{proposition}
{\bf Proof.} By Proposition $\ref{p-optimaltime}$, for $0\leq
s\leq t\leq T$
\begin{eqnarray*}
&&E_s^G[\sup_{\nu\in\mathcal{F}_{t,T}}E^G_t[f(S_{\nu})]]\\
&=&\sup_{\nu\in\mathcal{F}_{t,T}}E^G_s[f(S_{\nu})].
\end{eqnarray*}
Since $\mathcal{F}_{t,T}\subseteq \mathcal{F}_{s,T}$, we have
\begin{eqnarray*}
&&\sup_{\nu\in\mathcal{F}_{t,T}}E^G_s[f(S_{\nu})]\\
&\leq &\sup_{\nu\in\mathcal{F}_{s,T}}E^G_s[f(S_{\nu})].
\end{eqnarray*}
Thus, we derive that
\begin{eqnarray*}
&&E_s^G[\sup_{\nu\in\mathcal{F}_{t,T}}E^G_t[f(S_{\nu})]]\\
&\leq &\sup_{\nu\in\mathcal{F}_{s,T}}E^G_s[f(S_{\nu})].
\end{eqnarray*}
We prove the Proposition.$\ \ \ \ \ \ \ \ \square$

\begin{theorem} \label{th-bid-ask}Assume that the uncertain financial market consists of the bond which has the price process satisfying $(\ref{bond})$ and $d-$ risky assets with the
price processes as the G-asset price systems $(\ref{G-asset})$ and
can trade freely, the contingent claim $\xi$ which is written on
the $d$ assets with the maturity $T>0$ has the class of the payoff
defined in Definition $\ref{payoff}$, and the function
$V^{Am}(t,S_t)$ is defined in $(\ref{valuefunction})$. Then there
exists a superhedging  (resp. subhedging) strategy for $\xi$, such
that, the process $V=(V_t)_{0\leq t\leq T}$ defined by
\begin{eqnarray}
V_t:=e^{-r(T-t)}V^{Am}(t,S_t),& (\mbox{resp.  }
-e^{-r(T-t)}\mbox{ess
sup}_{\nu\in\mathcal{F}_{t,T}}E^G_t[-f_{\nu}])
\end{eqnarray}
is the ask (resp. bid) price process against $\xi$.
\end{theorem}
{\bf Proof.} The value function for the optimal stop time
$V^{Am}(t,S_t)$ is a G-supermartingale, it is easily to check that
$e^{-rt}V_t$ is G-supermartingale. By G-Doob-Meyer decomposition
Theorem $\ref{GDoob}$
\begin{eqnarray}
e^{-rt}V_t=M_t-\bar{C}_t
\end{eqnarray}
where $M_t$ is a G-martingale and $\bar{C}_t$ is an increasing
process with $\bar{C}_0=0$. By G-martingale representation Theorem
(\cite{Peng2010} and \cite{Song})
\begin{eqnarray}
M_t=E^G[M_T]+\int_0^t\eta_sd\tilde{B}_t-K_t
\end{eqnarray}
where $\eta_s\in H_G^1(0,T)$, $-K_t$ is a G-martingale, and $K_t$
is an increasing process with $K_0=0$. From the above equation, we
have
\begin{eqnarray}
e^{-rt}V_t= E^G[M_T]+\int_0^t\eta_sd\tilde{B}_t-(K_t+\bar{C}_t),
\end{eqnarray}
hence $(V_t,e^{rt}\eta_t,\int_0^te^{rs}d(\bar{C}_s+K_s)ds)$ is a
superhedging strategy.

Assume that $(Y_t,\pi_t, C_t)$ is a
superhedging strategy against $\xi$, then
\begin{eqnarray}\label{strategy2}
e^{-rt}Y_t=e^{-rT}\xi-\int_t^T\pi_td\tilde{B}_t+C_t.
\end{eqnarray}
Taking conditional G-expectation on the both sides of the equation
$(\ref{strategy2})$ and notice that the process $C_t$ is an
increasing process with $C_0=0$, we derive
\begin{eqnarray}
e^{-rt}Y_t&\geq &E_t^G[e^{-rT}\xi]
\end{eqnarray}
which implies that
\begin{eqnarray*}
Y_t&\geq &E_t^G[e^{-r(T-t)}\xi]\\
&\geq &E_t^G[e^{-r(T-t)}\mbox{ess sup}_{\nu\in
\mathcal{F}_{T,T}}[f_{\nu}]\\
&\geq & e^{-r(T-t)}\mbox{ess sup}_{\nu\in
\mathcal{F}_{T,T}}E_t^G[f_{\nu}]\\
&\geq & e^{-r(T-t)}\mbox{ess sup}_{\nu\in
\mathcal{F}_{t,T}}E_t^G[f_{\nu}]\\
&=& V_t
\end{eqnarray*}
from which, we prove that $V_t=e^{-r(T-t)}V^{Am}(t,S_t) $ is the
ask price against the claim $\xi$ at time $t$. Similarly we can
prove that $-e^{-r(T-t)}\mbox{ess sup}_{\nu\in
\mathcal{F}_{t,T}}E_t^G[-f_{\nu}]$ is the bid price against the
claim $\xi$ at time $t$. $\ \ \square$

\section{Free Boundary and Optimal Stopping Problems}

For given $t\in [0,T]$, $x\in R^d$, the G-asset price system
($\ref{G-asset}$) of the risky asset can be rewritten as follows
\begin{eqnarray*}
\left\{\begin{array}{l}
dS_{u}^{t,x}=S_{u}^{t,x}(rdt+d\tilde{B}_t)\\
S_{t}^{t,x}=x
\end{array}\right.
\end{eqnarray*}

We define the following deterministic function
\begin{eqnarray*}
u^a(t,x) &:=& e^{-r(T-t)}V^{Am}(t,S_t^{t,x}),
\end{eqnarray*}
where
\begin{eqnarray*}
V^{Am}(t,S_t^{t,x})=\mbox{ess sup}_{\nu\in
\mathcal{F}_{t,T}}E^G_{t} [f(S_{\nu}^{t,x})].
\end{eqnarray*}
From Theorem $\ref{th-bid-ask}$ the price of an American option
with expiry date $T$ and payoff function $f$, is the value
function of the optimal stopping problem
\begin{eqnarray}\label{stoppingproblem}
u^a(t,x) &:=& e^{-r(T-t)}\mbox{ess sup}_{\nu\in
\mathcal{F}_{t,T}}E^G_{t} [f(S_{\nu})].
\end{eqnarray}
We define operator $L$ as follows:
\begin{eqnarray*}
Lu=G(D^2u)+rD u+\partial_tu,
\end{eqnarray*}
where $G(\cdot)$ is the sublinear function defined by equation
$(\ref{sublinear-function})$. We consider the free boundary
problem
\begin{eqnarray}\label{freeboundary}\left\{
\begin{array}{ll}
\mathcal{L}u:=\max\{Lu-ru,f-u\}=0,& \mbox{in }[0,T]\times R^d,\\
u(T,\cdot)=f(T,\cdot),& \mbox{in } R^d
\end{array}\right.
\end{eqnarray}
Denote
$$
\mathcal{S}_T:=[0,T]\times R^d,
$$
for $p\geq 1$
$$
\mathcal{S}^p(\mathcal{S}_T):=\{u\in L^p(\mathcal{S}_T): D^2u,Du,
\partial_tu\in L^p(\mathcal{S}_T) \}.
$$
and for any compact subset $D$ of $\mathcal{S}_T$, we denote
$\mathcal{S}_{\mbox{loc}}^p(D)$ as the space of functions $u\in
\mathcal{S}^p(D)$.

\begin{definition} A function $u\in\mathcal{S}_{\mbox{loc}}^1(\mathcal{S}_T)\cap C(R^d\times
[0,T])$is a strong solution of problem $(\ref{freeboundary})$ if
$\mathcal{L}u=0$ almost everywhere in $\mathcal{S}_T$ and it
attains the final datum pointwisely. A function
$u\in\mathcal{S}_{\mbox{loc}}^1(\mathcal{S}_T)\cap C(R^d\times
[0,T])$ is a strong super-solution of problem
$(\ref{freeboundary})$ if $\mathcal{L}u\leq 0$.
\end{definition}
We will prove the following existence results
\begin{theorem} \label{freeboundaryexistence} If there exists a
strong super-solution $\bar{u}$ of problem $(\ref{freeboundary})$
then there also exists a strong solution $u$ of
$(\ref{freeboundary})$ such that $u\leq \bar{u}$ in
$\mathcal{S}_T$. Moreover $u\in
\mathcal{S}_{\mbox{loc}}^p(\mathcal{S}_T)$ for any $p\geq 1$ and
consequently, by the embedding theorem we have $u\in
C_{\mbox{B,loc}}^{1,\alpha}(\mathcal{S}_T)$ for any $\alpha\in
[0,1]$.
\end{theorem}
\begin{theorem}\label{representation}
Let $u$ be a strong solution to the free boundary problem
$(\ref{freeboundary})$ such that
\begin{eqnarray}\label{th-asumme-5-1}
|u(t,x)|\leq Ce^{\lambda|x|^2},&(t,x)\in \mathcal{S}_T,
\end{eqnarray}
form some constants $C,\lambda$ with $\lambda$ sufficiently small
so that
$$
E^G[\exp(\lambda\sup_{t\leq u\leq T}|S_u^{t,x}|^2)]<\infty
$$
holds. Then we have
\begin{eqnarray}
u(t,x)&=& e^{-r(T-t)}\mbox{ess sup}_{\nu\in
\mathcal{F}_{t,T}}E^G_{t} [f(S_{\nu})],
\end{eqnarray}
i.e., the solution of the free boundary problem is the value
function of the optimal stopping problem. In particular such a
solution is unique.
\end{theorem}

\subsection{Proof of Theorem $\ref{representation}$}

We employ a truncation and regularization technique to exploit the
weak interior regularity properties of $u$, for $R>0$ we set for
$R>0$, $B_R:=\{x\in R^d||x|<R \}$, and for $x\in B_R$ denote by
$\tau_R$ the first exit time of $S_u^{t,x}$ from $B_R$, it is easy
check that $E^G[\tau_R]$ is finite. As a first step we prove the
following result: for every $(t,x)\in [0,T]\times B_R$ and
$\tau\in \mathcal{F}_{t,T}$ such that $\tau\in [t,\tau_R]$, it
holds
\begin{eqnarray}\label{th-5-1}
u(t,x)=E^G[u(\tau,S_{\tau}^{t,x})]-E^G[\int_t^{\tau}Lu
(s,S_s^{t,x})ds].
\end{eqnarray}
For fixed, positive and small enough $\varepsilon$, we consider a
function $u^{\varepsilon, R}$ on $R^{d+1}$ with compact support
and such that $u^{\varepsilon, R}=u$ on $[t, T-\varepsilon]\times
B_R$. Moreover we denote by $(u^{\varepsilon, R,n})_{n\in N}$ a
regularizing sequence obtained by convolution of $u^{\varepsilon,
R}$ with the usual mollifiers, then for any $p\geq 1$ we have
$u^{\varepsilon, R,n}\in \mathcal{S}^p(R^{d+1})$ and
\begin{eqnarray}\label{1ststep}
\lim_{n\longrightarrow
\infty}\|Lu^{\varepsilon,R,n}-Lu^{\varepsilon,R}\|_{L^p([t,
T-\varepsilon]\times B_R)}=0.
\end{eqnarray}
By G-It$\hat{o}$ formula we have
\begin{eqnarray}
\begin{array}{rcl}
u^{\varepsilon,R,n} (\tau,S_{\tau}^{t,x})&=&u^{\varepsilon,R,n}
(t,x)+\frac{1}{2}\int_t^{\tau}D^2u^{\varepsilon,R,n}d<B>_s\\
&&+\int_t^{\tau}rDu^{\varepsilon,R,n}ds
+\int_t^{\tau}\partial_su^{\varepsilon,R,n}ds+\int_t^{\tau}Du^{\varepsilon,R,n}dB_s,
\end{array}
\end{eqnarray}
which implies that
\begin{eqnarray}
E^G[u^{\varepsilon,R,n} (\tau,S_{\tau}^{t,x})]=u^{\varepsilon,R,n}
(t,x)+\int_t^{\tau}Lu^{\varepsilon,R,n}ds.
\end{eqnarray}
We have
$$
\lim_{n\longrightarrow\infty}u^{\varepsilon,R,n}(t,x)=u^{\varepsilon,R}(t,x)
$$
and, by dominated convergence
$$
\lim_{n\longrightarrow\infty}E^G[u^{\varepsilon,R,n}(\tau,S_{\tau}^{t,x})]=E^G[u^{\varepsilon,R}(\tau,S_{\tau}^{t,x})].
$$
We have
\begin{eqnarray*}
&&|E^G[\int_t^{\tau}Lu^{\varepsilon,R,n}(s,S_{s}^{t,x})ds]-E^G[\int_t^{\tau}Lu^{\varepsilon,R}(s,S_{s}^{t,x})ds]|\\
&\leq &
E^G[\int_t^{\tau}|Lu^{\varepsilon,R,n}(s,S_{s}^{t,x})-Lu^{\varepsilon,R}(s,S_{s}^{t,x})|ds],
\end{eqnarray*}
by sublinear expectation representation Theorem (see
\cite{Peng2010}) there exists a family of probability space $Q$,
such that
\begin{eqnarray*}
&&E^G[\int_t^{\tau}|Lu^{\varepsilon,R,n}(s,S_{s}^{t,x})-Lu^{\varepsilon,R}(s,S_{s}^{t,x})|ds]\\
 &=& \mbox{ess sup}_{P\in Q}
E_P[\int_t^{\tau}|Lu^{\varepsilon,R,n}(s,S_{s}^{t,x})-Lu^{\varepsilon,R}(s,S_{s}^{t,x})|ds].
\end{eqnarray*}
Since $\tau\leq \tau_R$
\begin{eqnarray*}
&&\mbox{ess sup}_{P\in Q}
E_P[\int_t^{\tau}|Lu^{\varepsilon,R,n}(s,S_{s}^{t,x})-Lu^{\varepsilon,R}(s,S_{s}^{t,x})|ds]\\
&\leq & \mbox{ess sup}_{P\in Q}
E_P[\int_t^{T-\varepsilon}|Lu^{\varepsilon,R,n}(s,y)-Lu^{\varepsilon,R}(s,y)|I_{|S_{s}^{t,x}|\leq B_R}ds]\\
&\leq & \mbox{ess sup}_{P\in Q}
\int_t^{T-\varepsilon}\int_{B_R}|Lu^{\varepsilon,R,n}(s,y)-Lu^{\varepsilon,R}(s,y)|\Gamma_P(t,x;s,y)dyds
\end{eqnarray*}
where $\Gamma_P(t,x;\cdot,\cdot)\in L^{\bar{q}}([t,T]\times B_R)$,
for some $\bar{q}>1$, is the transition density of the solution of
$$
dX_s^{t,x}=X_s^{t,x}(rds+\sigma_{s,P}dW_{s,P})
$$
where $W_{s,P}$ is Wiener process in probability space
$(\Omega_t,P,\mathcal{F}^{P},\mathcal{F}_t^{P})$, and
$\sigma_{s,P}$ is adapted process such that $\sigma_{s,P}\in
[\underline{\sigma},\overline{\sigma}]$. By H$\ddot{o}$lder
inequality, we have ($1/\bar{p}+1/\bar{q}=1$)
\begin{eqnarray*}
\int_t^{T-\varepsilon}\int_{B_R}|Lu^{\varepsilon,R,n}(s,y)-Lu^{\varepsilon,R}(s,y)|\Gamma_P(t,x;s,y)dyds\leq
\|Lu^{\varepsilon,R,n}(s,y)-Lu^{\varepsilon,R}(s,y)\|_{L^{\bar{q}}([t,T]\times
B_R)}\|\Gamma_P(t,x;s,y)\|_{L^{\bar{p}}([t,T]\times B_R)},
\end{eqnarray*}
then, we obtain that
\begin{eqnarray*}
\lim_{n\longrightarrow \infty}E^G
[\int_t^{\tau}Lu^{\varepsilon,R,n}(s,S_{s}^{t,x})]=E^G[\int_t^{\tau}Lu^{\varepsilon,R}(s,S_{s}^{t,x})].
\end{eqnarray*}
This concludes the proof of $(\ref{th-5-1})$, since
$u^{\varepsilon,R}=u$ on $[t,T-\varepsilon]\times B_R$ and
$\varepsilon>0$ is arbitrary.

Since $Lu\leq 0$, we have for any $\tau\in \mathcal{F}_{t,T}$
$$
E^G\int_{t}^{\tau}Lu(s, S_s^{t,x})ds\leq 0.
$$
we infer from $(\ref{th-5-1})$ that
$$
u(t,x)\ge E^G[u(\tau\wedge\tau_R,S_{\tau\wedge \tau_R}^{t,x})].
$$
Next we pass to the limit as $R\longrightarrow +\infty$: we have
$$
\lim_{R\longrightarrow +\infty}\tau\wedge\tau_R=\tau,
$$
and by the growth assumption $(\ref{th-asumme-5-1})$
$$
|u(\tau\wedge\tau_R,S_{\tau\wedge \tau_R}^{t,x})|\leq
C\exp(\lambda \sup_{t\leq s\leq T}|S_{s}^{t,x}|^2).
$$
As $R\longrightarrow +\infty$
$$
u(t,x)\geq E^G[u(\tau, S_{\tau}^{t,x})]\geq E^G[f(\tau,
S_{\tau}^{t,x})].
$$
This shows that
$$
u(t,x)\geq \sup_{\tau\in \mathcal{F}_{t,T}}E^G[f(\tau,
S_{\tau}^{t,x})].
$$
We conclude the proof by putting
$$
\tau_0=\inf\{s\in[t,T]|u(s, S_{s}^{t,x})=f(s, S_{s}^{t,x})\}.
$$
Since $Lu=0$ a.e. where $u>\phi$, it holds
$$
E^G[\int_t^{\tau_0\wedge\tau_R}Lu(s,S_s^{t,x})ds]=0
$$
and from $(\ref{th-5-1})$ we derive that
$$
u(t,x)=E^G[u(\tau_0\wedge \tau_R, S_{\tau_0\wedge \tau_R}^{t,x})]
$$
Repeating the previous argument to pass to the limit in $R$, we
obtain
$$
u(t,x)=E^G[u(\tau_0,S_{\tau_0}^{t,x})]=E^G[f(\tau_0,S_{\tau_0}^{t,x})].
$$
Therefore, we finish the proof.   $\ \ \ \ \ \square$
\subsection{Free boundary problem}
Here we consider the free boundary problem on a bounded cylinder.
We denote the bounded cylinders as the form $[0,T]\times H_n$,
where $(H_n)$ is an increasing covering of $R^d$. We will prove
the existence of a strong solution to problem
\begin{eqnarray}\label{freeboundarycylinder}
\left\{\begin{array}{ll} \max\{Lu,f-u\}=0,&\mbox{in
}H(T):=[0,T]\times H,\\
u|_{\partial_PH(T)}=f,&
\end{array}
\right.
\end{eqnarray}
where $H$ is a bounded domain of $R^d$ and
\begin{eqnarray*}
\partial_PH(T):=\partial H(T)\setminus (\{T\}\times H)
\end{eqnarray*}
is the parabolic boundary of $H(T)$.

We assume the following condition on the payoff function
\begin{assum}\label{assum-payoff}
The payoff function $\xi=f(S_u^{t,x})$ have the following
assumption expressed by the sublinear function
\begin{eqnarray}
-G(-D^2 f)\geq c& \mbox{in } H,
\end{eqnarray}
where $G(\cdot)$ is the sublinear function defined by equation
$(\ref{sublinear-function})$.
\end{assum}
\begin{theorem}\label{th-freeboundarycylinder}
We assume the Assumption $\ref{assum}$ holds. Problem
$(\ref{freeboundarycylinder})$ has a strong solution $u\in
\mathcal{S}_{loc}^1(H(T))\cap C(\overline{H(T)})$. Moreover $u\in
\mathcal{S}_{loc}^p(H(T))$ for any $p>1$.
\end{theorem}
{\bf Proof.} The proof is based on a standard penalization
technique (see Friedeman \cite{Friedman}). We consider a family
$(\beta_{\varepsilon})_{\varepsilon\in [0,1]}$ of smooth functions
such that, for any $\varepsilon$, the function
$\beta_{\varepsilon}$ is increasing, bounded on $R$ and has
bounded first order derivative, such that
\begin{eqnarray*}
\beta_{\varepsilon}(s)\leq \varepsilon,\ \ s>0, \ \ \mbox{and}\ \
\lim_{\varepsilon\longrightarrow
0}\beta_{\varepsilon}(s)=-\infty,\ \ s<0.
\end{eqnarray*}
We denote by $f^{\delta}$ as the regularization of $f$, and
consider the following penalized  and regularized problem and
denote the solution as $u_{\varepsilon,\delta}$
\begin{eqnarray}\label{regularfreeboundary}
\left\{\begin{array}{ll}
Lu=\beta_{\varepsilon}(u-f^{\delta}),&\mbox{in
}H(T),\\
u|_{\partial_PH(T)}=f^{\delta},&
\end{array}
\right.
\end{eqnarray}
Lions \cite{Lions}, Krylov \cite{Krylov} and Nisio \cite{Nisio}
prove that problem ($\ref{regularfreeboundary}$) has a unique
viscosity solution $u_{(\varepsilon,\delta)}\in
C^{2,\alpha}(\overline{H(T)})\cap C(\overline{H(T)})$ with
$\alpha\in [0,1]$.

Next, we firstly prove the uniform boundedness of the penalization
term:
\begin{eqnarray}\label{penalizationterm}
|\beta_{\varepsilon}(u_{\varepsilon,\delta}-f^{\delta})|\leq
c,&\mbox{ in } H(T),
\end{eqnarray}
with $c$ independent of $\varepsilon$ and $\delta$.

By construction $\beta_{\varepsilon}\leq \varepsilon$, it suffices
to prove the lower bound in ($\ref{penalizationterm}$). By
continuity,
$\beta_{\varepsilon}(u_{\varepsilon,\delta}-f^{\delta})$ has a
minimum $\zeta$ in $\overline{H(T)}$ and we may suppose
$$
\beta_{\varepsilon}(u_{\varepsilon,\delta}(\zeta)-f^{\delta}(\zeta))\leq
0,
$$
otherwise we prove the lower bound. If $\zeta\in \partial_{P}H(T)$
then
$$
\beta_{\varepsilon}(u_{\varepsilon,\delta}(\zeta)-f^{\delta}(\zeta))=\beta_{\varepsilon}(0)=
0.
$$
On the other hand, if $\zeta\in H(T)$, then we recall that
$\beta_{\varepsilon}$ is increasing and consequently
$u_{(\varepsilon,\delta)}-f^{\delta}$ also has a (negative)
minimum in $\zeta$. Thus, we have
\begin{eqnarray}\label{bound}
Lu_{\varepsilon,\delta}(\zeta)-Lf^{\delta}(\zeta)\geq 0\geq
u_{\varepsilon,\delta}(\zeta)-f^{\delta}(\zeta).
\end{eqnarray}
By the Assumption 5.1 on $f$, we have that $Lf^{\delta}(\zeta)$ is
bounded uniformly in $\delta$. Therefor, by $(\ref{bound})$, we
deduce
\begin{eqnarray*}
\beta_{\varepsilon}(u_{(\varepsilon,\delta)}(\zeta)-f^{\delta}(\zeta))&=&Lu_{(\varepsilon,\delta)}(\zeta)\\
&\geq& Lf^{\delta}(\zeta)\geq c,
\end{eqnarray*}
where $c$ is a constant independent on $\varepsilon,\delta$ and
this proves $(\ref{penalizationterm})$.

Secondly, we use the $\mathcal{S}^p$ interior estimate combined
with $(\ref{penalizationterm})$, to infer that, for every compact
subset $D$ in $H(T)$ and $p\geq 1$, the norm
$\|u_{\varepsilon,\delta}\|_{\mathcal{S}^p(D)}$ is bounded
uniformly in $\varepsilon$ and $\delta$. It follows that
$(u_{\varepsilon,\delta})$ converges as
$\varepsilon,\delta\longrightarrow 0$, weakly in $\mathcal{S}^p$
on compact subsets of $H(T)$ to a function $u$. Moreover
$$
\limsup_{\varepsilon,\delta}\beta_{\varepsilon}(u_{\varepsilon,\delta}-f^{\delta})\leq
0,
$$
so that $Lu\leq f$ a.e. in $H(T)$. On the other hand, $Lu=f$ a.e.
in the set $\{u>f\}$.

Finally, it is straightforward to verify that $u\in
C(\overline{H(T)})$ and assumes the initial-boundary conditions,
by using standard arguments based on the maximum principle and
barrier functions.$\ \ \ \ \ \square$

{\bf Proof of Theorem $\ref{freeboundaryexistence}$.}The proof of
Theorem $\ref{freeboundaryexistence}$ about the existence theorem
for the free boundary problem on unbounded domains is similar in
\cite{Francesco} by using Theorem $\ref{th-freeboundarycylinder}$
about the existence theorem for the free boundary problem on the
regular bounded cylindrical domain.$\ \ \ \ \ \square$

\end{document}